\magnification=\magstep1          % |
\hsize=16truecm                   % | don't change these values !!
\vsize=23.5truecm                 % |
%

%  set \Finalout at the top of your paper, if you don't like
%  overfullrules!
\parindent=6mm
\normallineskiplimit=.99pt
\catcode`\@=11
%
%fonts (preloaded in PLAIN and AMS fonts)
%roman:
 \font\ninrm=cmr9 \font\egtrm=cmr8 \font\sixrm=cmr6
%math italic
 \font\nini=cmmi9 \skewchar\nini='177
 \font\egti=cmmi8 \skewchar\egti='177
 \font\sixi=cmmi6 \skewchar\sixi='177
%symbols
 \font\ninsy=cmsy9 \skewchar\ninsy='60
 \font\egtsy=cmsy8 \skewchar\egtsy='60
 \font\sixsy=cmsy6 \skewchar\sixsy='60
 \font\ninex=cmex9 \font\egtex=cmex8
%bold
 \font\ninbf=cmbx9 \font\egtbf=cmbx8 \font\sixbf=cmbx6
%text italic
 \font\ninit=cmti9 \font\egtit=cmti8
%slanted
 \font\ninsl=cmsl9 \font\egtsl=cmsl8
%caps and small caps
 \font\tensc=cmcsc10 \font\ninsc=cmcsc9 \font\egtsc=cmcsc8
%sans serif
 \font\tensf=cmss10 \font\ninsf=cmss9 \font\egtsf=cmss8
%euler and blackboard-bold
% \font\teneum=eufm10 \font\seveneum=eufm7 \font\fiveeum=eufm5   % EUL10
% \font\nineum=eufm9 \font\sixeum=eufm6                          % EUL9
% \font\tenmsb=msbm10 \font\sevenmsb=msbm7 \font\fivemsb=msbm5   % BBB10
% \font\ninmsb=msbm9 \font\sixmsb=msbm6                          % BBB9
%title
 \font\titbf=cmbx10 scaled \magstep2
 \font\titsyt=cmbsy10 scaled \magstep2 \skewchar\titsyt='60
 \font\titsys=cmbsy7 scaled \magstep2 \skewchar\titsys='60
 \font\titsyss=cmbsy5 scaled \magstep2 \skewchar\titsyss='60
 \font\titmit=cmmib10 scaled \magstep2 \skewchar\titmit='177
 \font\titmis=cmmib7 scaled \magstep2 \skewchar\titmis='177
 \font\titmiss=cmmib5 scaled \magstep2 \skewchar\titmiss='177
 \font\titit=cmbxti10 scaled \magstep2
% \font\titeut=eufb10 scaled \magstep2                           % EUL14
% \font\titeus=eufb7 scaled \magstep2                            % EUL14
% \font\titeuss=eufb5 scaled \magstep2                           % EUL14
% \font\titmsbt=msbm10 scaled \magstep2                          % BBB14
% \font\titmsbs=msbm7 scaled \magstep2                           % BBB14
% \font\titmsbss=msbm5 scaled \magstep2                          % BBB14
%sectionheadings
 \font\secbf=cmbx10 scaled \magstep1
 \font\secsyt=cmbsy10 scaled \magstep1 \skewchar\secsyt='60
 \font\secsys=cmbsy7 scaled \magstep1 \skewchar\secsys='60
 \font\secsyss=cmbsy5 scaled \magstep1 \skewchar\secsyss='60
 \font\secmit=cmmib10 scaled \magstep1 \skewchar\secmit='177
 \font\secmis=cmmib7 scaled \magstep1 \skewchar\secmis='177
 \font\secmiss=cmmib5 scaled \magstep1 \skewchar\secmiss='177
 \font\secit=cmbxti10 scaled \magstep1
% \font\seceut=eufb10 scaled \magstep1                           % EUL12
% \font\seceus=eufb7 scaled \magstep1                            % EUL12
% \font\seceuss=eufb5 scaled \magstep1                           % EUL12
% \font\secmsbt=msmb10 scaled \magstep1                          % BBB12
% \font\secmsbs=msbm7 scaled \magstep1                           % BBB12
% \font\secmsbss=msbm5 scaled \magstep1                          % BBB12
%subsectionheadings
 \font\ssecsyt=cmbsy10 \skewchar\ssecsyt='60
 \font\ssecsys=cmbsy7 \skewchar\ssecsys='60
 \font\ssecsyss=cmbsy5 \skewchar\ssecsyss='60
 \font\ssecmit=cmmib10 \skewchar\ssecmit='177
 \font\ssecmis=cmmib7 \skewchar\ssecmis='177
 \font\ssecmiss=cmmib5 \skewchar\ssecmiss='177
 \font\ssecit=cmbxti10
% \font\sseceut=eufb10 \font\sseceus=eufb7 \font\sseceuss=eufb5  % EUL10
%
%new families: euler and blackboard bold
%\newfam\eufam                                                   % EUL
%\def\frak{%                                                     % .
% \let\next\relax                                                % .
% \ifmmode
%  \let\next\frak@
% \else
%  \def\next{\errmessage{Use \string\frak\space in math mode only }}
% \fi
% \next}
%\def\frak@#1{{\frak@@{#1}}}                                     % .
%\def\frak@@#1{\fam\eufmfam#1}                                   % .
%\let\goth\frak                                                  % EUL
%
%\newfam\bbfam                                                   % BBB
%\def\bbb{%                                                      % .
% \let\next\relax                                                % .
% \ifmmode\let\next\bbb@
% \else\def\next{\errmessage{Use \string\bbb\space in math mode only}}
% \fi
% \next}                                                         % .
%\def\bbb@#1{{\bbb@@{#1}}}                                       % .
%\def\bbb@@#1{\fam\bbfam#1}                                      % BBB
%
%font sizes (adapted from MANMAC.TEX)
\newdimen\b@gsize
\def\b@g#1#2{{\hbox{$\left#2\vcenter to#1\b@gsize{}%
                    \right.\n@space$}}}
\def\big{\b@g\@ne}
\def\Big{\b@g{1.5}}
\def\bigg{\b@g\tw@}
\def\Bigg{\b@g{2.5}}
\def\tenpoint{%
 \textfont0=\tenrm \scriptfont0=\sevenrm \scriptscriptfont0=\fiverm%
 \def\rm{\fam0\tenrm}%
 \textfont1=\teni \scriptfont1=\seveni \scriptscriptfont1=\fivei%
 \textfont2=\tensy \scriptfont2=\sevensy \scriptscriptfont2=\fivesy%
 \textfont3=\tenex \scriptfont3=\tenex \scriptscriptfont3=\tenex%
 \def\it{\fam\itfam\tenit}\textfont\itfam=\tenit%
 \def\sl{\fam\slfam\tensl}\textfont\slfam=\tensl%
 \textfont\bffam=\tenbf \scriptfont\bffam=\sevenbf%
 \scriptscriptfont\bffam=\fivebf%
 \def\bf{\fam\bffam\tenbf}%
 \def\sc{\tensc}%
 \def\sf{\tensf}%
% \textfont\eufam=\teneum \scriptfont\eufam=\seveneum%           % EUL10
% \scriptscriptfont\eufam=\fiveeum%                              % EUL10
% \textfont\bbfam=\tenmsb \scriptfont\bbfam=\sevenmsb%           % BBB10
% \scriptscriptfont\bbfam=\fivemsb%                              % BBB10
 \normalbaselineskip=12pt%
 \setbox\strutbox=\hbox{\vrule height8.5pt depth3.5pt width\z@}%
 \abovedisplayskip12pt plus3pt minus9pt%
 \belowdisplayskip12pt plus3pt minus9pt%
 \abovedisplayshortskip\z@ plus3pt%
 \belowdisplayshortskip7pt plus3pt minus4pt%
 \normalbaselines\rm%
 \setbox\z@\vbox{\hbox{$($}\kern\z@}\b@gsize=1.2\ht\z@%
}
\def\ninepoint{%
 \textfont0=\ninrm \scriptfont0=\sixrm \scriptscriptfont0=\fiverm%
 \def\rm{\fam0\ninrm}%
 \textfont1=\nini \scriptfont1=\sixi \scriptscriptfont1=\fivei%
 \textfont2=\ninsy \scriptfont2=\sixsy \scriptscriptfont2=\fivesy%
 \textfont3=\ninex \scriptfont3=\ninex \scriptscriptfont3=\ninex%
 \def\it{\fam\itfam\ninit}\textfont\itfam=\ninit%
 \def\sl{\fam\slfam\ninsl}\textfont\slfam=\ninsl%
 \textfont\bffam=\ninbf \scriptfont\bffam=\sixbf%
 \scriptscriptfont\bffam=\fivebf%
 \def\bf{\fam\bffam\ninbf}%
 \def\sc{\ninsc}%
 \def\sf{\ninsf}%
% \textfont\eufam=\nineum \scriptfont\eufam=\sixeum%             % EUL9
% \scriptscriptfont\eufam=\fiveeum%                              % EUL9
% \textfont\bbfam=\ninmsb \scriptfont\bbfam=\sixmsb%             % BBB9
% \scriptscriptfont\bbfam=\fivemsb%                              % BBB9
 \normalbaselineskip=11pt%
 \setbox\strutbox=\hbox{\vrule height8pt depth3pt width\z@}%
 \abovedisplayskip11pt plus2.7pt minus8.1pt%
 \belowdisplayskip11pt plus2.7pt minus8.1pt%
 \abovedisplayshortskip\z@ plus2.7pt%
 \belowdisplayshortskip6.3pt plus2.7pt minus3.6pt%
 \normalbaselines\rm%
 \setbox\z@\vbox{\hbox{$($}\kern\z@}\b@gsize=1.2\ht\z@%
}
\def\@eightpoint{%
 \textfont0=\egtrm \scriptfont0=\sixrm \scriptscriptfont0=\fiverm%
 \def\rm{\fam0\egtrm}%
 \textfont1=\egti \scriptfont1=\sixi \scriptscriptfont1=\fivei%
 \textfont2=\egtsy \scriptfont2=\sixsy \scriptscriptfont2=\fivesy%
 \textfont3=\egtex \scriptfont3=\egtex \scriptscriptfont3=\egtex%
 \def\it{\fam\itfam\egtit}\textfont\itfam=\egtit%
 \def\sl{\fam\slfam\egtsl}\textfont\slfam=\egtsl%
 \textfont\bffam=\egtbf \scriptfont\bffam=\sixbf%
 \scriptscriptfont\bffam=\fivebf%
 \def\bf{\fam\bffam\egtbf}%
 \def\sc{\egtsc}%
 \def\sf{\egtsf}%
 \normalbaselineskip=9pt%
 \setbox\strutbox=\hbox{\vrule height7pt depth2pt width\z@}%
 \abovedisplayskip10pt plus2.4pt minus7.2pt%
 \belowdisplayskip10pt plus2.4pt minus7.2pt%
 \abovedisplayshortskip\z@ plus2.4pt%
 \belowdisplayshortskip5.6pt plus2.4pt minus3.2pt%
 \normalbaselines\rm%
 \setbox\z@\vbox{\hbox{$($}\kern\z@}\b@gsize=1.2\ht\z@%
}
\def\@titlefonts{%
 \textfont0=\titbf \def\rm{\titbf}%
 \textfont1=\titmit%
 \scriptfont1=\titmis \scriptscriptfont1=\titmiss%
 \textfont2=\titsyt%
 \scriptfont2=\titsys \scriptscriptfont2=\titsyss%
 \textfont\itfam=\titit \def\it{\titit}%
% \textfont\eufam=\titeut \scriptfont\eufam=\titeus%             % EUL14
% \scriptscriptfont\eufam=\titeuss%                              % EUL14
% \textfont\bbfam=\titmsbt \scriptfont\bbfam=\titmsbs%           % BBB14
% \scriptscriptfont\bbfam=\titmsbss%                             % BBB14
 \let\bf\rm \let\sf\rm \let\sc\rm \let\sl\it%
 \normalbaselineskip18pt%
 \setbox\strutbox\hbox{%
  \vrule height 12.6pt depth5.4pt width\z@}%
 \normalbaselines\rm%
 \setbox\z@\vbox{\hbox{$($}\kern\z@}\b@gsize=1.2\ht\z@%
}
\def\@sectionfonts{%
 \textfont0=\secbf \def\rm{\secbf}%
 \textfont1=\secmit%
 \scriptfont1=\secmis \scriptscriptfont1=\secmiss%
 \textfont2=\secsyt%
 \scriptfont2=\secsys \scriptscriptfont2=\secsyss%
 \textfont\itfam=\secit \def\it{\secit}%
% \textfont\eufam=\seceut \scriptfont\eufam=\seceus%             % EUL12
% \scriptscriptfont\eufam=\seceuss%                              % EUL12
% \textfont\bbfam=\secmsbt \scriptfont\bbfam=\secmsbs%           % BBB12
% \scriptscriptfont\bbfam=\secmsbss%                             % BBB12
 \let\bf\rm \let\sf\rm \let\sc\rm \let\sl\it%
 \normalbaselineskip14.5pt%
 \setbox\strutbox\hbox{%
  \vrule height 10.1pt depth4.4pt width\z@}%
 \normalbaselines\rm%
 \setbox\z@\vbox{\hbox{$($}\kern\z@}\b@gsize=1.2\ht\z@%
}
\def\@subsecfonts{\tenpoint%
 \textfont0=\tenbf \def\rm{\tenbf}%
 \textfont1=\ssecmit%
 \scriptfont1=\ssecmis \scriptscriptfont1=\ssecmiss%
 \textfont2=\ssecsyt%
 \scriptfont2=\ssecsys \scriptscriptfont2=\ssecsyss%
 \textfont\itfam=\ssecit \def\it{\ssecit}%
% \textfont\eufam=\sseceut \scriptfont\eufam=\sseceus%           % EUL10
% \scriptscriptfont\eufam=\sseceuss%                             % EUL10
 \let\bf\rm \let\sf\rm \let\sc\rm \let\sl\it%
 \normalbaselines\rm%
}
\def\vfootnote#1{%
  \insert\footins\bgroup%
  \interlinepenalty\interfootnotelinepenalty%
  \ninepoint%
  \splittopskip\ht\strutbox \splitmaxdepth\dp\strutbox%
  \floatingpenalty\@MM%
  \leftskip\z@skip \rightskip\z@skip \spaceskip\z@skip \xspaceskip\z@skip%
  \textindent{#1}\footstrut\futurelet\next\fo@t}
\def\makefootline{\baselineskip18pt\line{\the\footline}}
\newif\ifSec
\outer\def\section#1\par{\@beginsection{#1}}
\def\@beginsection#1{%
 \tenpoint%
 \vskip\z@ plus.1\vsize\penalty-250%
 \vskip\z@ plus-.1\vsize\vskip1.75\bigskipamount%
 \message{#1}\leftline{\noindent\@sectionfonts#1}%
 \nobreak\bigskip\noindent%
 \global\Sectrue
 \everypar{\ifSec\global\Secfalse\everypar{}\fi}}%
\outer\def\references{\@beginsection{References}\parindent20pt\ninepoint}%
\let\ref\item%
\def\subsection#1{\tenpoint
 \ifSec\else\bigskip\fi
 {\@subsecfonts#1}}%
\long\def\proclaim#1{%
 \medbreak
 \bf#1\enspace\it\ignorespaces}
\long\def\statement#1{%
 \medbreak
 \bf#1\enspace\rm\ignorespaces}
\def\endproclaim{%
 \ifdim\lastskip<\medskipamount
  \removelastskip\penalty55\medskip\fi
 \rm}
\let\endstatement\endproclaim
\def\qedrule{\hbox{\vrule height1.4ex depth0pt width1ex}}
\def\qed{\unskip\nobreak\quad\qedrule\medbreak}
\newskip\figskipamount \figskipamount\bigskipamount
\newbox\@caption
\def\figskip{\vskip\figskipamount}
\long\def\figure#1#2{%
 \setbox\@caption\hbox{\ninepoint\rm\ignorespaces#2}%
 \ifhmode%
  \vadjust{%
   \figskip%
   \line{\hfill\vbox to #1{\vfill}}%
   \ifvoid\@caption\else%
    \vskip\medskipamount%
    \ifdim \wd\@caption>\hsize%
     \noindent{\ninepoint\rm\ignorespaces#2\par}%
     \setbox\@caption\hbox{}%
    \else%
     \line{\hfil\box\@caption\hfil}%
    \fi%
   \fi%
   \figskip}%
 \else \ifvmode
  \figskip
  \line{\hfill\vbox to #1{\vfill}}
  \ifvoid\@caption\else
   \medskip
   \ifdim \wd\@caption>\hsize
    \noindent{\ninepoint\rm\ignorespaces#2\par}
    \setbox\@caption\hbox{}
   \else
    \line{\hfil\box\@caption\hfil}
   \fi
  \fi
  \figskip
 \fi\fi}
\newcount\@startpage
\@startpage=1
\def\@logo{\vtop{\@eightpoint
 \line{\strut\hfill}
 \line{\strut\hfill}
 \line{\strut\hfill}}}
\def\logo#1#2#3#4#5{\global\@startpage=#4 \global\pageno=#4
 \def\@logo{\vtop{\@eightpoint
 \line{\hfill Zeitschrift f\"ur Analysis und ihre Anwendungen}
 \line{\hfill Journal for Analysis and its Applications}
 \line{\hfill Volume #1 (#2), No.\ #3, #4--#5}}}}
\headline={%
 \ifnum\pageno=\@startpage
  \hfill
 \else
  \ifodd\pageno
   \strut\hfill{\ninepoint \sf \@runtitle}\qquad{\tenpoint\bf\folio}
  \else
   \strut{\tenpoint\bf\folio}\qquad{\ninepoint\sf\@runauthor}\hfill
  \fi
 \fi}
\footline={%
 \ifnum\pageno=\@startpage
  \hfill
 \else
  \hfill
 \fi}
\def\author#1{\def\@author{\ignorespaces#1}}
\def\@author{}
\def\runauthor#1{\def\@runauthor{\ignorespaces#1\ et.~al.}}
\def\@runauthor{\@author}
\def\title#1{\def\@title{\def\\{\hfill\egroup\line\bgroup\strut\hfill}#1}}
\def\@title{}
\def\runtitle#1{\def\@runtitle{\ignorespaces#1}}
\def\@runtitle{\@title}
\long\def\abstract#1\endabstract{\def\@abstract{\ignorespaces#1}}
\def\@abstract{}
\def\keywords#1{\def\@keywords{\ignorespaces#1}}
\def\@keywords{}
\def\primclass#1{%
 \def\@pclass{\ignorespaces#1}
 \global\let\classification\undefined
 \gdef\@class{}}
\def\secclasses#1{\def\@sclass{\ignorespaces#1}}
\def\classification#1{%
 \def\@class{\ignorespaces#1}
 \global\let\primclass\undefined
 \global\let\secclasses\undefined
 \gdef\@pclass{}\gdef\@sclass{}}
\newcount\addresscount
\addresscount=1
\newcount\addressnum
\def\address#1{%
 \expandafter\gdef\csname @address\number\addresscount \endcsname {#1}
 \global\advance\addresscount by 1}
\def\@addresses{
 \kern-10pt
 \addressnum=0
 \loop
  \ifnum\addressnum<\addresscount
   \advance\addressnum by 1
   \csname @address\number\addressnum \endcsname
   \ifnum\addressnum<\addresscount\hfill\break\fi
 \repeat
}
\long\def\maketitle{\begingroup
 \parindent0pt
\def\afootnote{%
  \insert\footins\bgroup%
  \interlinepenalty\interfootnotelinepenalty%
  \ninepoint%
  \splittopskip\ht\strutbox \splitmaxdepth\dp\strutbox%
  \floatingpenalty\@MM%
  \leftskip\z@skip \rightskip\z@skip \xspaceskip\z@skip%
  \footstrut\futurelet\next\fo@t}
 \vglue-12.5mm
 \@logo
 \vskip16.5pt
 \begingroup\@titlefonts
 \line\bgroup\strut\hfill\@title\hfill\egroup
 \endgroup
 \bigskip
 \centerline{\tenpoint\noindent\bf\hfill\@author\hfill}
 \vskip20mm
 \vbox{\ninepoint\noindent \bf Abstract.\ \rm\@abstract\par}
 \medskip
 \vbox{\ninepoint\noindent \bf Keywords:\ \it\@keywords\par}
 \smallskip
 \vbox{\ninepoint\noindent \bf AMS subject classification:\
  \ifx\classification\undefined
   \rm Primary\ \@pclass, secondary\ \@sclass\par
  \else
   \rm \@class
  \fi}
 \afootnote{\@addresses}
 \let\maketitle\relax
 \gdef\@title{} \gdef\@abstract{}
 \gdef\@keywords{}
 \gdef\@pclass{} \gdef\@sclass{} \gdef\@class{}
 \addressnum=0
 \loop
  \ifnum\addressnum<\addresscount
  \advance\addressnum by 1
  \expandafter\gdef\csname @address\number\addressnum \endcsname{}
 \repeat
 \endgroup
}
% switch to tenpoint
\catcode`@=12
\tenpoint
%%%%%%%%%%%%%%%%%%%%%%%%%%%%%%%%%%%%%%%%%%%%%%%%%%%%%
%  end of macro   / beginning of TeX-file
%%%%%%%%%%%%%%%%%%%%%%%%%%%%%%%%%%%%%%%%%%%%%%%%%%%%%
\author{M.~Frank}
\address{M.~Frank: University of Leipzig, Inst.~Math., D-04109 Leipzig,
F.~R.~G.}
\title{A multiplier approach to the Lance-Blecher theorem}

\runtitle{On the Lance-Blecher theorem}

\abstract
A new approach to the Lance-Blecher theorem is presented resting on the
interpretation of elements of Hilbert C*-module theory in terms of multiplier
theory of operator C*-algebras:
The Hilbert norm on a Hilbert C*-module allows to recover the values of the
inducing C*-valued inner product in a unique way, and two Hilbert C*-modules
$\{ {\cal M}_1, \langle .,. \rangle_1 \}$, $\{ {\cal M}_2, \langle .,.
\rangle_2 \}$ are isometrically iso\-mor\-phic as Banach C*-modules if and
only if there exists a bijective C*-linear map $S: {\cal M}_1 \to {\cal M}_2$
such that the identity $\langle .,. \rangle_1 \equiv \langle S(.),S(.)
\rangle_2$ is valid. In particular, the values
of a C*-valued inner product on a Hilbert C*-module are completely determined
by the Hilbert norm induced from it.
In addition, we obtain that two C*-valued inner products
on a Banach C*-module inducing equivalent norms to the given one give rise
to isometrically isomorphic Hilbert C*-modules if and only if the derived
C*-algebras of ''compact'' module operators are $*$-isomorphic. The
involution and the C*-norm of the C*-algebra of ''compact'' module operators
on a Hilbert C*-module allow to recover its original C*-valued inner product
up to the following equivalence relation:
$\langle .,. \rangle_1 \sim \langle .,. \rangle_2$ if and only if there
exists an invertible, positive element $\, a \,$ of the center of ${\rm M}(A)$
such that the identity $\langle .,. \rangle_1 \equiv a \cdot \langle .,.
\rangle_2$
holds.
\endabstract


\primclass{46L99}
\secclasses{46H25}

\maketitle

%\vfill \eject

%%% BEGINNING OF TEXT %%%
\smallskip \noindent
One of the crucial problems of the theory of Hilbert C*-modules is the
interrelation of three kinds of isomorphisms between them: Banach C*-module
isomorphisms, isometric Banach C*-module isomorphisms and unitary C*-module
isomorphisms intertwining the C*-valued inner products. Already in 1985
L.~G.~Brown gave examples of pairs of C*-valued inner products on a
certain Banach C*-module which induce equivalent norms, but without
any linking isometric Banach C*-module automorphism between these Hilbert
C*-modules, [3, Ex.~6.2, 6.3] (cf.~[13, Ex.~2.3], [4]).
In 1994 E.~C.~Lance showed that two Hilbert C*-modules $\{ {\cal M}_1, \langle
.,. \rangle_1 \}$, $\{ {\cal M}_2, \langle .,. \rangle_2 \}$ are isometrically
isomorphic as Banach C*-modules if and only if there exists a bijective
C*-linear map $S: {\cal M}_1 \to {\cal M}_2$ such that the identity
$\langle .,. \rangle_1 \equiv \langle S(.),S(.) \rangle_2$ is valid on
${\cal M}_1 \times {\cal M}_1$, [11, Th.]. The theorem of E.~C.~Lance
indicates that the Hilbert norm on a Hilbert C*-module might determine the
possible values of the related C*-valued inner product(s), in the best case
up to uniqueness.
Looking for examples one recalls that for every Hilbert space the concrete
values of its inner product can be exactly recovered from the Hilbert norm.
Surprisingly, the same turns out to be true for general Hilbert C*-modules
as D.~P.~Blecher obtained in 1995, [1, Th.~3.1, 3.2], [2]. He used the
point of view of operator modules over (non-self-adjoint) operator algebras
and elements of the representation theories of C*-algebras and of Hilbert
C*-modules to prove this result.
\smallskip \noindent
We want to give an alternative purely C*-algebraic proof of this important
fact pointing out the related background in multiplier theory of C*-algebras.
The formula (1) telling how to recover the values of the C*-valued inner
product from the Hilbert norm was partially suggested by the approach of
D.~P.~Blecher. Beside this, our approach has the advantage that intertwining
isomorphisms of C*-valued inner products on a Banach C*-module which induce
equivalent norms to the given one can be expressed in terms of the related
operator C*-algebras: Two Hilbert C*-modules $\{ {\cal M}, \langle .,.
\rangle_1
\}$, $\{ {\cal M}, \langle .,. \rangle_2 \}$ are isometrically isomorphic as
Banach C*-modules if and only if the derived C*-algebras of ''compact''
module operators are $*$-isomorphic. For this aim results on quasi-multipliers
of
C*-algebras are involved due to L.~G.~Brown, [3]. The involution and the
C*-norm of the C*-algebra of ''compact'' module operators on a Hilbert
C*-module allow to recover its original C*-valued inner product up to
the following equivalence relation:
$\langle .,. \rangle_1 \sim \langle .,. \rangle_2$ if and only if there
exists an invertible, positive element $\, a \,$ of the center of ${\rm M}(A)$
such that the identity $\langle .,. \rangle_1 \equiv a \cdot \langle .,.
\rangle_2$
holds for arbitrary elements of the given Hilbert C*-module.
If the center of ${\rm M}(A)$ is trivial then one has only to fix the
Hilbert norm on one singular non-zero element of the Hilbert C*-module to make
the choice unique.
%-----------------------------------------------------------------------------

\smallskip \noindent
We start our investigations recalling some definitions and basic facts
from the literature, cf.~[8,9,10,12,15].
\noindent
We consider Hilbert C*-modules $\{ {\cal M}, \langle .,. \rangle \}$ over
general C*-algebras $A$, i.~e.~(left) $A$-modules $\cal M$ together with
an $A$-valued inner product $\langle .,. \rangle: {\cal M} \times {\cal M}
\rightarrow A$ satisfying the conditions:
\smallskip
(i) $\langle x,x \rangle \ge 0$ for every $x \in {\cal M}$.
\smallskip
(ii) $\langle x,x \rangle =0$ if and only if $x=0$.
\smallskip
(iii) $\langle x,y \rangle = \langle y,x \rangle^*$ for every $x,y \in
{\cal M}$.
\smallskip
(iv) $\langle ax+by,z \rangle = a\langle x,z \rangle + b \langle y,z \rangle$
for
every $a,b \in A$, $x,y,z \in {\cal M}$.
\smallskip
(v) ${\cal M}$ is complete with respect to the norm $\|x\|=\|\langle x,x
\rangle \|_A^{1/2}$.
\smallskip \noindent
We always suppose, that the linear structures of the C*-algebra $A$ and of
the (left) $A$-module ${\cal M}$ are compatible, i.~e. $\lambda(ax)=(\lambda
a)x=
a(\lambda x)$ for every $\lambda \in {\bf C}$, $a \in A$, $x \in {\cal M}$.
A Hilbert C*-module is said to be full if the norm-closed linear span of
the values of the C*-valued inner product coincides with its C*-algebra
of coefficients.
Let us denote the $A$-dual Banach $A$-module of a Hilbert $A$-module
$\{ {\cal M}, \langle .,. \rangle \}$ by ${\cal M}'$, where ${\cal M}'=
\{ r:{\cal M} \to A \, : \, r \; {\rm is} \; \hbox{\rm {\it A}-linear} \;
{\rm and} \; {\rm bounded} \}$. A Hilbert C*-module $\{ {\cal M}, \langle .,.
\rangle \}$ is self-dual if the standard isometric C*-linear embedding $x \in
{\cal M} \to \langle .,x \rangle \in {\cal M}'$ is surjective.
\noindent
The class of (self-dual) Hilbert W*-modules is of special interest. Many
pathologies can be avoided for them because the C*-valued inner product lifts
always to the C*-dual Banach W*-module turning it into a self-dual Hilbert
W*-module, [15]. To each Hilbert C*-module $\{ {\cal M},\langle .,.
\rangle \}$ over a C*-algebra $A$ one can assign a standard Hilbert W*-module
over the bidual W*-algebra $A^{**}$ of $A$ in the following way,
cf.~[14, Def.~1.3], [15, \S 4]:
Form the algebraic tensor product $A^{**}\otimes {\cal M}$ which becomes a
(left) $A^{**}$-module defining the action of $A^{**}$ on its
elementary tensors by the formula $ab \otimes x = a (b \otimes x)$ for
$a,b \in A^{**}$, $x \in \cal M$. Now , setting
$$
\left[ \sum_i a_i \otimes x_i , \sum_j b_j \otimes y_j \right] =
\sum_{i,j} a_i \langle x_i, y_j \rangle b_j
$$
on finite sums of elementary tensors one obtains a degenerate
$A^{**}$-valued inner pre-product. The factorization of $A^{**}
\otimes {\cal M}$ by the set $\{ z \in A^{**} \otimes {\cal M} : [z,z]=0 \}$
gives a Hilbert $A^{**}$-module denoted by ${\cal M}^{\#}$ in the sequel.
It contains $\cal M$ as a $A$-submodule.
If $\cal M$ is self-dual then ${\cal M}^{\#}$ is self-dual, too, but the
converse conclusion is still an open problem.
Every bounded $A$-linear operator $T$ on $\cal M$ has a unique extension
to a bounded $A^{**}$-linear operator on ${\cal M}^{\#}$ preserving
the operator norm.

\smallskip \noindent
In the following we want to consider several kinds of module operators on
Hilbert C*-modules. An $A$-linear bounded operator
$K$ on a Hilbert $A$-module $\{ {\cal M}, \langle .,. \rangle \}$ is
''compact'' if it belongs to the norm-closed linear hull of the set of
elementary operators
$\{ \theta_{x,y} \, : \, \theta_{x,y}(z) = \langle z,x \rangle y \, , \,
x,y \in {\cal M} \} \, ,$
[10,15]. The set of all ''compact'' operators on ${\cal M}$ is denoted by
${\rm K}_A({\cal M})$. A bounded $A$-linear operator on a Hilbert C*-module
${\cal M}$ is adjointable if the operator $T^*$ defined by the formula
$\langle T(x),y \rangle = \langle x,T^*(y) \rangle$ for all $x,y \in {\cal M}$
is a bounded $A$-linear operator on ${\cal M}$.
By [9,10] the C*-algebra ${\rm K}_A({\cal M})$ is a two-sided ideal of the set
of all
bounded, adjointable module operators ${\rm End}_A^*({\cal M})$ on ${\cal M}$
which is $*$-isomorphic to its multiplier C*-algebra.
\smallskip \noindent
To characterize unitary isomorphisms of Hilbert C*-modules we use the
following definition:

\statement{Definition 1.}
 \noindent
 Let $A$ be a fixed C*-algebra. Two Hilbert $A$-modules $\{ {\cal M}_1,
 \langle .,. \rangle_1 \}$ and $\{ {\cal M}_2, \langle .,. \rangle_2 \}$ are
 said to be {\it isomorphic as Hilbert C*-modules} (or equiv., {\it unitarily
 isomorphic}) if and only if there exists a linear bijective mapping
 $S: {\cal M}_1 \rightarrow {\cal M}_2$ such that the equalities
 $S(ax) = a S(x)$ and $\langle x,y \rangle_1 = \langle S(x),S(y) \rangle_2$
 are valid for every $a \in A$, every $x,y \in {\cal M}_1$.
\endstatement

\noindent
The literature contains some results about the existence of such isomorphisms
between Hilbert C*-modules:
If a Hilbert $A$-module $\{ {\cal M},\langle .,. \rangle_1 \}$ over a
given C*-algebra $A$ is self-dual then every $A$-valued inner
product $\langle .,. \rangle _2$ on $\cal M$ inducing an equivalent
to the given one norm fulfills the identity $\langle .,. \rangle_2 = \langle
S(\cdot),S(\cdot) \rangle_1$ on ${\cal M} \times \cal M$ for a unique
positive invertible bounded $A$-linear operator $S$ on $\cal M$,
cf.~[5, Th.~2.6]. Similarly, E.~C.~Lance proved for arbitrary Hilbert
$A$-modules ${\cal M}_1$, ${\cal M}_2$ over a fixed C*-algebra $A$ that in
case of existence of a bounded $A$-linear adjointable operator $T: {\cal M}_1
\to {\cal M}_2$ with dense ranges for $T$ and $T^*$ there exists an unitary
Banach C*-module isomorphism of ${\cal M}_1$ and ${\cal M}_2$, [12, Prop.~3.8].
Countably generated Hilbert C*-modules are isomorphic as Banach C*-modules
if and only if they are isometrically isomorphic as Banach C*-modules,
[3, Cor.~4.8, Th.~4.9], [6, Th.~3.1].

\noindent
To explain what kind of general results one could obtain we prefer
to rely on multiplier theory of C*-algebras. The fundamental result
of H.~Lin cited below appears to be very helpful. (In fact, it extends a
well-known result of P.~Green and G.~G.~Kasparov.)

\proclaim{Proposition 2.} {\rm ([13, Th.~1.5, 1.6], cf.~[9,10] and [17]) }

     \noindent
     Let $A$ be a C*-algebra and let $\{ {\cal M} , \langle .,. \rangle \}$
     be a Hilbert $A$-module. Then the mapping $\phi$ defined by the
     formula
     $$
     \phi : {\rm End}_A({\cal M},{\cal M}') \rightarrow {\rm QM}({\rm K}_A
     ({\cal M})) \, ,
     \, \theta_{x,y} \phi(T) \theta_{z,t} = \theta_{(T(t)(x))z,y}
     $$
     \noindent
     $(x,y,z,t \in {\cal M})$ is an isometric isomorphism of involutive Banach
     spaces.

     \noindent
     The restriction of $\phi$ to ${\rm End}_A(\cal M)$ induces an isometric
     algebraic isomorphism to the Banach algebra ${\rm LM}({\rm K}_A(\cal M))$.

     \noindent
     Finally, the restriction of $\phi$ to ${\rm End}_A^*(\cal M)$ induces a
     $*$-isomorphism to the C*-algebra ${\rm M}({\rm K}_A(\cal M))$.
\endproclaim

\noindent
Note, that every left Hilbert $A$-module $\cal M$ can be considered as a
right Hilbert ${\rm K}_A({\cal M})$-module fixing another ${\rm K}_A
({\cal M})$-valued inner product $\langle x,y \rangle_{Op.} = \theta_{x,y}$.
This point of view gives another interpretation of the left actions of
${\rm M}(A)$ and of ${\rm LM}(A)$ on full Hilbert $A$-modules $\cal M$ by
Proposition 2.

\smallskip \noindent
Our key observation is that every $A$-valued inner product $\langle .,.
\rangle$ on a Hilbert C*-module $\cal M$ over a given C*-algebra $A$ defines a
mapping $T$ from $\cal M$ into its $A$-dual Banach $A$-module ${\cal M}'$
by the formula $T: x \in {\cal M} \rightarrow \langle .,x \rangle \in
{\cal M}'$. The properties of these mappings $T$ in terms of multiplier
C*-theory are the following:

\proclaim{Proposition 3.}
     Let $A$ be a C*-algebra and let $\{ {\cal M} , \langle .,. \rangle_1
     \}$ be a Hilbert $A$-module. Denote by $\langle .,. \rangle_2$ a
     second $A$-valued inner product on $\cal M$ inducing an equivalent
     norm to the given one. Then the mapping $T: x \in {\cal M} \rightarrow
     \langle .,x \rangle_2 \in {\cal M}'$ can be identified with a uniquely
     defined invertible positive element of ${\rm QM}({\rm K}_A({\cal M}))
     \subset {\rm K}_A^{(1)}({\cal M})^{**}$. Conversely, every invertible
     positive element $T' \in {\rm QM}({\rm K}_A({\cal M})) \subset
     {\rm K}_A^{(1)}({\cal M})^{**}$ induces an $A$-valued inner product
     and an equivalent norm on $\cal M$ via the formula $\langle x,y \rangle_2
     = (\phi^{-1}(T')(y))(x)$ for $x,y \in \cal M$.
\endproclaim

\statement{Proof.}  Using the identifications of Proposition 2 made by the
mapping $\phi$ one derives the equality
$$
\theta^{(1)}_{x,y} \phi(T) \theta^{(1)}_{z,t} = \theta^{(1)}_{\langle x,t
\rangle_2 z,y} \in {\rm K}^{(1)}_A({\cal M}) \, ,
$$
\noindent
which defines $\phi(T) \in {\rm QM}({\rm K}_A({\cal M}))$ by the right side
of this equality.
To show the positivity of the quasi-multiplier $\phi(T)$ one modifies the
equality above setting $x=t$, $y=z$. Making use of the identity $\theta_{x,y}
=\theta_{y,x}^*$ valid for every $x,y \in \cal M$ one obtains
$$
\langle \theta^{(1)}_{\langle x,x \rangle_2 t,t}(s), s \rangle_1
\, = \, \langle \langle s,t \rangle_1 \langle z,z \rangle_2 t,s \rangle_1
\, = \,  \langle s,t \rangle_1 \langle z,z \rangle_2 \langle t,s \rangle_1
\, \geq \, 0
$$
\noindent
for every $s \in \cal M$. Since $\phi(T) \in {\rm K}^{(1)}_A(\cal M)^{**}$ by
construction and since the linear span of the ''compact'' operators of type
$\theta$ is
norm-dense inside ${\rm K}^{(1)}_A(\cal M)$ the positivity of $\phi(T)$ as an
element of the W*-algebra ${\rm K}^{(1)}_A(\cal M)^{**}$ follows.

\smallskip \noindent
To show the invertibility of $\phi(T)$ inside ${\rm K}^{(1)}_A(\cal M)^{**}$
we use a standard construction from the introduction. First,
build the Hilbert $A^{**}$-module ${\cal M}^{\#}$ from $\cal M$.
Both the $A$-valued inner products on $\cal M$ can be extended in a unique
way to $A^{**}$-valued inner products on ${\cal M}^{\#}$. Secondly, take
the (self-dual) $A^{**}$-dual Hilbert $A^{**}$-module $({\cal M}^{\#})'$ of
${\cal M}^{\#}$. Again, both the inner products can be continued
(cf.~[15, Th.~3.2, 3.6]), and their extensions are connected by an invertible
positive operator $S$ as described in [5, Th.~2.6]. Obviously,
the uniquely defined extension of the operator $T$ inside
${\rm End}_{A^{**}}(({\cal M}^{\#})')$ equals $S^*S$. Hence, the
(real) spectrum of $T$ is deleted away from zero by a positive constant, and
$\phi(T)$ is invertible.

\smallskip \noindent
Conversely, set $\langle x,y \rangle_2 = (\phi^{-1}(T')(y))(x)$ for $x,y \in
\cal M$ and for a given invertible positive $T' \in {\rm QM}({\rm K}_A
({\cal M})) \subset {\rm K}_A^{(1)}({\cal M})^{**}$. As can be easily
seen by consi\-derations similar to that above $\langle .,. \rangle _2$ is an
$A$-valued inner product on $\cal M$ inducing an equivalent norm to the
given one.  \qed
\endstatement
%%%%%%%%%%%%%%%%%%%%%%%%%%%%%%%bearbeitungsende
\statement{Example 4.}
     Let $A$ be a C*-algebra. Define the action of $A$ on itself by
     multiplication from the left. Then $A$ becomes a Hilbert $A$-module
     setting $\langle a,b \rangle_T = aTb^*$ for every $a,b \in A$ and a
     fixed positive invertible $T \in {\rm QM}(A)$. Vice versa, every
     $A$-valued inner product on $A$ arises in this manner.
     If $A$ is unital then $T \in A \equiv {\rm QM}(A)$.
\endstatement

\proclaim{Theorem 5.}  {\rm (E.~C.~Lance [11, Th.] /
                                      D.~P.~Blecher, [1, Th.~3.1,3.2])}

 \noindent
 Let $A$ be a C*-algebra and $\cal M$ be a left Banach $A$-module the norm of
 which is known to be generated by an $A$-valued inner product on $\cal M$
 with unknown values. Then this $A$-valued inner product $\langle .,. \rangle$
 on $\cal M$ is unique, and it can be recovered by the formulae
 $$
 \langle x,x \rangle := \sup \{ r(x)r(x)^* : r \in {\cal M}', \|r\| \leq 1 \}
 \, , \eqno(1)
 $$
 $$
 \langle x,y \rangle := {1 \over 4} \sum_{k=0}^3 {\rm i}^k \langle x+
 {\rm i}^ky, x+{\rm i}^ky \rangle
 $$
 \noindent
 for every $x,y \in \cal M$, where the right side of (1) uses the norm
 of the underlying Banach $A$-module only.

 \noindent
 Consequently, every bijective isometric $A$-linear isomorphism of two
 Hilbert $A$-modules $S: \{ {\cal M}_1, \langle .,. \rangle_1 \} \to
 \{ {\cal M}_2, \langle .,. \rangle_2 \}$ identifies the two $A$-valued
 inner products by the formula $\langle .,. \rangle_1 \equiv
 \langle S(.),S(.) \rangle_2$, and vice versa.
\endproclaim

\statement{Proof.} Let $\langle .,. \rangle_1$, $\langle .,. \rangle_2$ be two
$A$-valued inner products on $\cal M$ giving the norm. Applying again the
standard construction from the introduction both the $A$-valued
inner products can be continued to $A^{**}$-valued inner products on the
self-dual Hilbert $A^{**}$-module $({\cal M}^\#)'$. Inside ${\rm End}_{A^{**}}
(({\cal M}^\#)')$ exists a positive invertible operator $T$ such that the
identity $\langle .,. \rangle_2 \equiv \langle T(\cdot),. \rangle_1$ holds for
the continued $A^{**}$-valued inner products on $({\cal M}^\#)' \times
({\cal M}^\#)'$, cf.~[5, Prop.~2.2]. By construction one has
$$
\|x\| \equiv \|\langle T(x),x \rangle_1\|_A
                                   \equiv \|\langle x,x \rangle_1\|_A \quad ,
\quad \|x\| \equiv \|\langle T^{-1}(x),x \rangle_2\|_A
                                   \equiv \|\langle x,x \rangle_2\|_A \, ,
$$
\noindent
and by a theorem of W.~L.~Paschke [15, Th.~2.8] one obtains
$$
\|\langle T(x),x \rangle_1\|_A
              \leq \|T^{1/2}\|^2 \cdot \|\langle x,x \rangle_1\|_A \quad ,
              \quad
\|\langle T^{-1}(x),x \rangle_2\|_A
              \leq \|T^{-1/2}\|^2 \cdot \|\langle x,x \rangle_2\|_A \, .
$$
\noindent
This implies $\|T\|=\|T^{-1}\|=1$, and by the positivity of $T$ and by
general spectral properties of elements of C*-algebras $T={\rm id}_{\cal M}$
yields.

\smallskip \noindent
To show the estimation formula (1) of the values of the
$A$-valued inner product one recalls that
$$
r(x)r(x)^* \leq \|r\|^2 \langle x,x \rangle
$$
\noindent
for every $x \in {\cal M}$, every $r \in {\cal M}'$ by [15, Th.~2.8].
Since the $A$-valued inner product $\langle .,. \rangle$ induces an
isometric $A$-linear embedding of $\cal M$ into ${\cal M}'$ by the
formula $x \rightarrow \langle .,x \rangle$ one has only to indicate a
sequence $\{ r_n : n \in {\bf N} \}$ of bounded by one $A$-linear
functionals on $\cal M$ of this special nature such that the expressions
$\{ r_n(x)r_n(x)^* : n \in {\bf N} \}$ converge
to the value $\langle x,x \rangle$ in norm from below. This can be done
setting $r_n= \langle .,(\langle x,x \rangle + 1/n \cdot 1_A)^{-1/2} x \rangle$
for $n \in {\bf N}$. Consequently, the supremum really exists. The second
formula is obvious.

\smallskip \noindent
To show the last statement one has to consider the two $A$-valued inner
products $\langle .,. \rangle_1$ and $\langle S(.),S(.) \rangle_2$ on the
Hilbert $A$-module ${\cal M}_1$ inducing exactly the same norm for every
element
of ${\cal M}_1$. Therefore, they take identically the same values by the
first
part of the proof.
\qed
\endstatement

\proclaim{Proposition 6.}
Let $A$ be a C*-algebra and $\cal M$ be a Banach $A$-module possessing two
$A$-valued inner products $\langle .,. \rangle_1$, $\langle .,. \rangle_2$
inducing equivalent to the given one norms on $\cal M$. Suppose,
$0 < C,D < \infty$ are the minimal real numbers for which the inequality
$\|x\|_1 \leq C \cdot \|x\|_2 \leq D \cdot \|x\|_1$ is satisfied for every
$x \in \cal M$. Then the inequality
$$
\langle x,x \rangle_1 \leq C \cdot \langle x,x \rangle_2
                          \leq D \cdot \langle x,x  \rangle_1 \eqno(2)
$$
\noindent
is valid for every $x \in \cal M$ and the same real numbers $C,D$.
\endproclaim

\statement{Proof.} We extend both the $A$-valued inner products to
$A^{**}$-valued
inner products to the self-dual Hilbert $A^{**}$-module $({\cal M}^\#)'$
using the standard construction. Then there exists a positive invertible
operator $T$ inside ${\rm End}_{A^{**}}(({\cal M}^\#)')$ such that the
identity $\langle .,. \rangle_2 \equiv \langle T(\cdot),. \rangle_1$ holds for
the continued $A^{**}$-valued inner products on $({\cal M}^\#)' \times
({\cal M}^\#)'$, cf.~[5].
Applying [15, Th.~2.8] to the operators $T$ and $T^{-1}$ on $({\cal M}^\#)'$
one obtains the minimal real numbers $C=\|T^{-1}\|$ and $D=\|T\|\cdot\|T^{-1}\|$
for which the inequality (2) is valid, and these constants
equal the minimal constants obtained in the comparison inequality of
the two norms. \qed
\endstatement

\statement{Remark 7.}
The expression (1) could make sense for more general Banach
C*-modu\-les than Hilbert C*-modules. However, if it would be well-defined for
every element $x$ of a Banach, non-Hilbert C*-module $\cal M$ then it should
be non-C*-linear and/or degenerated, anyway.
\endstatement

\noindent
For more similar results we refer the reader to the work of D.~P.~Blecher
who has treated Hilbert C*-modules as operator spaces and operator modules
over (non-self-adjoint) ope\-rator algebras using mainly geometric notions
like complete contractability and complete boundedness of mappings, for
example, [1,2]. The advantage of our approach comes to light in the
following statement characterizing isometric isomorphisms of different
C*-valued inner products on a fixed Banach C*-module in terms of
$*$-isomorphisms of the related operator C*-algebras.

\proclaim{Theorem 8.}
Let $A$ be a C*-algebra and let $\{ {\cal M}, \langle .,. \rangle_1\}$ be
a Hilbert $A$-module. Let $\langle .,. \rangle_2$ be another $A$-valued
inner product on $\cal M$ inducing an equivalent to the given one norm.
The following conditions are equivalent:
\smallskip
(i) The $A$-valued inner product $\langle .,. \rangle_2$ on $\cal M$
    is generated by an invertible bounded $A$-linear operator $S$ on
    ${\cal M}$ satisfying the identity $\langle .,. \rangle_2 \equiv
    \langle S(.), S(.) \rangle_1$ on ${\cal M} \times {\cal M}$.
\smallskip
(ii) The positive invertible quasi-multiplier $T$ of ${\rm K}^{(1)}_A
     ({\cal M})$ corresponding to the $A$-valued inner product $\langle .,.
     \rangle_2$ by Proposition 3 is decomposable as $T=S^*S$ for at least
     one invertible left multiplier $S$ of ${\rm K}^{(1)}_A({\cal M})$.
\smallskip
(iii) The C*-algebra ${\rm K}^{(2)}_A({\cal M})$ of ''compact''
      operators on ${\cal M}$ corresponding to the $A$-valued inner product
      $\langle .,. \rangle_2$ is $*$-isomorphic to the original C*-algebra
      of ''compact'' operators ${\rm K}^{(1)}_A({\cal M})$.
\smallskip
(iv) The C*-algebra ${\rm End}^{*,(2)}_A({\cal M})$ of
     adjointable bounded $A$-linear ope\-rators on ${\cal M}$ corresponding to
     the
     $A$-valued inner product $\langle .,. \rangle_2$ is $*$-isomorphic to
     the original C*-algebra of adjointable bounded $A$-linear operators
     ${\rm End}^{*,(1)}_A({\cal M})$.
\endproclaim

\statement{Proof.}
The implications (i)$\leftrightarrow$(ii) follow from Proposition 2 together
with the key Proposition 3. Keeping in mind Proposition 3 one adapts
L.~G.~Brown's results on quasi-multipliers [3, Th.~4.2, Prop.~4.4]
of (non-unital) C*-algebras to the C*-algebra ${\rm K}^{(1)}_A({\cal M})$:
For a positive invertible quasi-multiplier $T$ of ${\rm K}^{(1)}_A({\cal M})$
the C*-subalgebra $T^{1/2}{\rm K}^{(1)}_A({\cal M}) T^{1/2}$ of the bidual
W*-algebra ${\rm K}^{(1)}_A({\cal M})^{**}$ is $*$-isomorphic to
${\rm K}^{(1)}_A({\cal M})$ if and only if there exists a left multiplier $S$
of ${\rm K}^{(1)}_A({\cal M})$ such that $T=S^*S$ inside ${\rm K}^{(1)}_A
({\cal M})^{**}$. Thus, one obtains the equivalence of the conditions (ii) and
(iii), cf.~[4],[6]. The equivalence of the last two statements is shown in [7].
\qed
\endstatement

\statement{Example 9.}
Every positive invertible quasi-multiplier $T$ of a (non-unital) C*-algebra
$A$ is decomposable as $T=S^*S$ for at least one invertible left multiplier
$S$ of $A$ if and only if every pair of $A$-valued inner products
$\langle .,. \rangle_1$, $\langle .,. \rangle_2$ on $A$ inducing equivalent
norms to the given C*-norm is connected by an isometric Banach $A$-module
isomorphism $S$ of the two corresponding (left) Banach $A$-modules $\{ A,
\|.\|_1 \}$, $\{ A, \|.\|_2 \}$, cf.~[13, Ex.~2.3] for a counterexample.
\endstatement

\noindent
Note, that the equivalence of the conditions of Theorem 8 does not hold
any longer if one considers C*-valued inner products on different Banach
C*-modules and $*$-isomorphisms of corresponding operator C*-algebras,
in general. A counterexample can be found in [6,7].
The canonical question arising is whether the original Hilbert norm can
be recovered from the C*-norm of the related operator C*-algebras, or not.
The answer is given by the following statement.

\proclaim{Proposition 10.}
Let $A$ be a C*-algebra and $\{ {\cal M}, \langle .,. \rangle_1 \}$ be a
full Hilbert $A$-module possessing a second $A$-valued inner product
$\langle .,. \rangle_2$ which induces an equivalent norm to the given one.
Suppose, both the $A$-valued inner products define the same bounded
$A$-linear operators on ${\cal M}$ to be ''compact'', and both they
induce the same involution and C*-norm on this algebra of all ''compact''
$A$-linear operators. Then there exists an invertible positive element $\, a\,$
of the center of the multiplier C*-algebra ${\rm M}(A)$ of $A$ such that
the identity $\langle .,. \rangle_1 \equiv a \cdot \langle .,. \rangle_2$
holds on ${\cal M} \times {\cal M}$.

\noindent
If the center of ${\rm M}(A)$ is trivial then the condition $\|x\|=1$ for
some fixed non-zero $x \in {\cal M}$ makes the choice of the $A$-valued
inner product on ${\cal M}$ unique.
\endproclaim

\statement{Proof.}
Since both the related C*-algebras of ''compact'' operators coincide,
i.e. they are $*$-isomorphic, Theorem 8 applies: The invertible positive
quasi-multiplier $T$ corresponding to the $A$-valued inner product
$\langle .,. \rangle_2$ is decomposable as $T=S^*S$ for an invertible
left multiplier $S$ which can be considered as a bounded $A$-linear
operator on ${\cal M}$ by Proposition 3. In particular, the inequality
  $$ \langle K(x),x \rangle_2 = \langle (SK)(x),S(x) \rangle_1 \geq 0 $$
holds for every positive ''compact'' operator $K$, every $x \in {\cal M}$.
Consequently, $S$ commutes with every positive ''compact'' operator
and belongs to the center of the multiplier C*-algebra ${\rm End}_A^*
({\cal M})$
of ${\rm K}_A({\cal M})$. But, ${\rm Z}({\rm End}_A^*({\cal M}))$
consists of the operators $\{ a \cdot {\rm id}_{{\cal M}} : a \in
{\rm Z}({\rm M}(A)) \}$, and it is $*$-isomorphic to ${\rm Z}({\rm M}(A))$.
No further restrictions apply to $S$ and $T=S^*S$ since $\|x\|_1=\|a^{-1/2}
\cdot
x \|_2$ and $\|K(x)\|_1 = \|K(a^{-1/2} \cdot x)\|_2$ for every $x \in{\cal M}$,
every $K \in {\rm K}_A({\cal M})$ (where $T=a \in {\rm Z}({\rm M}(A))$).
\qed
\endstatement

\statement{Acknowledgement.}
The author is indebted to Nien-Tsu Shen for submitting a copy of her
Ph.D.~thesis to him. He is grateful to D.~P.~Blecher, L.~G.~Brown and
E.~C.~Lance for sending copies of their preprints and for comments on
the subject. The special net of C*-linear functionals in the proof of
Proposition 5 was obtained during discussions with A.~S.~Mishchenko.
A.~Kasparek brought the Hilbert space case of Proposition 10
to the authors attention.
\endstatement

\references

\item{[1]} Blecher,~D.~P.:
{\it A new approach to Hilbert C*-modules.}
preprint, Univ.~of Houston, Houston, Texas, U.S.A., 1995.
\smallskip

\item{[2]} Blecher,~D.~P.:
{\it A generalization of Hilbert modules.}
J.~Functional Analysis 136(1996), 365-421.
\smallskip

\item{[3]} Brown,~L.~G.:
{\it Close hereditary C*-subalgebras and the structure of quasi-multipliers.}
MSRI preprint \# 11211-85, Purdue University, West Lafayette, USA, 1985.
\smallskip

\item{[4]} Brown,~L.~G., J.~A.~Mingo and Nien-Tsu Shen:
{\it Quasi-multipliers and embeddings of Hilbert C*-bimodules.}
Canad.~J.~Math. 46(1994), 1150-1174.
\smallskip

\item{[5]} Frank,~M.:
{\it Self-duality and C*-reflexivity of Hilbert C*-modules.}
Zeitschr.~Anal.~Anwendungen 9(1990), 165-176.
\smallskip

\item{[6]} Frank,~M.:
{\it Geometrical aspects of Hilbert C*-modules.} preprint 22/93, K{\o}benhavns
Universitet, Matematisk Institut, 1993. submitted to Annals Global Anal.
Geom.in Sept.~1996.
\smallskip

\item{[7]} Frank,~M.:
{\it Isomorphisms of Hilbert C*-modules and $*$-isomorphisms of related
operator algebras.} accepted by Math.~Scand..
\smallskip

\item{[8]} Frank,~M.:
{\it Hilbert C*-modules and their applications -- a guided reference
overview.}
preprint 13/1996, ZHS-NTZ, Univ.~Leipzig, 1996.
\smallskip

\item{[9]} Green,~P.:
{\it The local structure of twisted covariance algebras.},
Acta Math. 140(1978), 191-250.
\smallskip

\item{[10]} Kasparov,~G.~G.:
{\it Hilbert $C^*$-modules: The theorems of Stinespring and Voiculescu.}
J.~Ope\-rator Theory 4(1980), 133-150.
\smallskip

\item{[11]} Lance,~E.~C.:
{\it Unitary operators on Hilbert C*-modules.}
Bull.~London Math.~Soc.~26(1994), 363-366.
\smallskip

\item{[12]} Lance,~E.~C.:
{\it Hilbert C*-modules -- a toolkit for operator algebraists.}
London Mathematical Society Lecture Note Series 210, Cambridge University
Press, Cambridge, England, 1995.
\smallskip

\item{[13]} Lin,~Huaxin:
{\it Bounded module maps and pure completely positive maps.}
J.~Operator Theory 26(1991), 121-138.
\smallskip

\item{[14]} Lin,~Huaxin:
{\it Injective Hilbert C*-modules.}
Pacific J.~Math. 154(1992), 131-164.
\smallskip

\item{[15]} Paschke,~W.~L.:
{\it Inner product modules over B*-algebras.} Trans.~Amer.~Math.~Soc.
182(1973), 443- 468.
\smallskip

\item{[16]} Shen, Nien-Tsu:
{\it Embeddings of Hilbert bimodules.} Ph.D., Purdue University, West
Lafa\-yette, U.S.A., 1982.

\item{[17]} Zeller-Meier,~G.:
{\it Some remarks about C*-Hilbert spaces and Hilbert C*-modules.}
preprint, Marseille, France, 1991.
\vskip10mm

Received on July 30, 1996.

\bye